\documentclass{article}
\usepackage{amssymb,amsmath, theorem,euscript,
 }

\newcounter{sec}

\def\SS{\smallskip}

\newcounter{punct}[sec]

\def\punct{\refstepcounter{punct}{\arabic{sec}.\arabic{punct}.  }}

\def\COUNTERS{\addtocounter{sec}{1}
              \setcounter{punct}{0}
          \setcounter{equation}{0}
          \setcounter{theorem}{0}
                  }

\newtheorem{theorem}{Theorem}[sec]
\newtheorem{proposition}[theorem]{Proposition}

\newtheorem{lemma}[theorem]{Lemma}

\newtheorem{corollary}[theorem]{Corollary}
\newtheorem{observation}[theorem]{Observation}

\begin{document}

\def\ov{\overline}
\def\wt{\widetilde}

\newcommand{\rk}{\mathop {\mathrm {rk}}\nolimits}
\newcommand{\Aut}{\mathop {\mathrm {Aut}}\nolimits}
\newcommand{\Out}{\mathop {\mathrm {Out}}\nolimits}
\renewcommand{\Re}{\mathop {\mathrm {Re}}\nolimits}

\def\Br{\mathrm {Br}}

\def\SL{\mathrm {SL}}
\def\SU{\mathrm {SU}}
\def\GL{\mathrm  {GL}}
\def\U{\mathrm  U}
\def\OO{\mathrm  O}
\def\Sp{\mathrm  {Sp}}
\def\SO{\mathrm  {SO}}
\def\SOS{\mathrm {SO}^*}
\def\Diff{\mathrm{Diff}}
\def\Vect{\mathfrak{Vect}}

\def\PGL{\mathrm  {PGL}}
\def\PU{\mathrm {PU}}

\def\PSL{\mathrm  {PSL}}

\def\Symp{\mathrm{Symp}}

\def\cA{\mathcal A}
\def\cB{\mathcal B}
\def\cC{\mathcal C}
\def\cD{\mathcal D}
\def\cE{\mathcal E}
\def\cF{\mathcal F}
\def\cG{\mathcal G}
\def\cH{\mathcal H}
\def\cJ{\mathcal J}
\def\cI{\mathcal I}
\def\cK{\mathcal K}
\def\cL{\mathcal L}
\def\cM{\mathcal M}
\def\cN{\mathcal N}
\def\cO{\mathcal O}
\def\cP{\mathcal P}
\def\cQ{\mathcal Q}
\def\cR{\mathcal R}
\def\cS{\mathcal S}
\def\cT{\mathcal T}
\def\cU{\mathcal U}
\def\cV{\mathcal V}
\def\cW{\mathcal W}
\def\cX{\mathcal X}
\def\cY{\mathcal Y}
\def\cZ{\mathcal Z}

 \def\cBS{\mathcal {B}}



\def\frA{\mathfrak A}
\def\frB{\mathfrak B}
\def\frC{\mathfrak C}
\def\frD{\mathfrak D}
\def\frE{\mathfrak E}
\def\frF{\mathfrak F}
\def\frG{\mathfrak G}
\def\frH{\mathfrak H}
\def\frJ{\mathfrak J}
\def\frK{\mathfrak K}
\def\frL{\mathfrak L}
\def\frM{\mathfrak M}
\def\frN{\mathfrak N}
\def\frO{\mathfrak O}
\def\frP{\mathfrak P}
\def\frQ{\mathfrak Q}
\def\frR{\mathfrak R}
\def\frS{\mathfrak S}
\def\frT{\mathfrak T}
\def\frU{\mathfrak U}
\def\frV{\mathfrak V}
\def\frW{\mathfrak W}
\def\frX{\mathfrak X}
\def\frY{\mathfrak Y}
\def\frZ{\mathfrak Z}

\def\fra{\mathfrak a}
\def\frb{\mathfrak b}
\def\frc{\mathfrak c}
\def\frd{\mathfrak d}
\def\fre{\mathfrak e}
\def\frf{\mathfrak f}
\def\frg{\mathfrak g}
\def\frh{\mathfrak h}
\def\fri{\mathfrak i}
\def\frj{\mathfrak j}
\def\frk{\mathfrak k}
\def\frl{\mathfrak l}
\def\frm{\mathfrak m}
\def\frn{\mathfrak n}
\def\fro{\mathfrak o}
\def\frp{\mathfrak p}
\def\frq{\mathfrak q}
\def\frr{\mathfrak r}
\def\frs{\mathfrak s}
\def\frt{\mathfrak t}
\def\fru{\mathfrak u}
\def\frv{\mathfrak v}
\def\frw{\mathfrak w}
\def\frx{\mathfrak x}
\def\fry{\mathfrak y}
\def\frz{\mathfrak z}

\def\frsp{\mathfrak{sp}}


\def\bfa{\mathbf a}
\def\bfb{\mathbf b}
\def\bfc{\mathbf c}
\def\bfd{\mathbf d}
\def\bfe{\mathbf e}
\def\bff{\mathbf f}
\def\bfg{\mathbf g}
\def\bfh{\mathbf h}
\def\bfi{\mathbf i}
\def\bfj{\mathbf j}
\def\bfk{\mathbf k}
\def\bfl{\mathbf l}
\def\bfm{\mathbf m}
\def\bfn{\mathbf n}
\def\bfo{\mathbf o}
\def\bfp{\mathbf p}
\def\bfq{\mathbf q}
\def\bfr{\mathbf r}
\def\bfs{\mathbf s}
\def\bft{\mathbf t}
\def\bfu{\mathbf u}
\def\bfv{\mathbf v}
\def\bfw{\mathbf w}
\def\bfx{\mathbf x}
\def\bfy{\mathbf y}
\def\bfz{\mathbf z}

\def\bfA{\mathbf A}
\def\bfB{\mathbf B}
\def\bfC{\mathbf C}
\def\bfD{\mathbf D}
\def\bfE{\mathbf E}
\def\bfF{\mathbf F}
\def\bfG{\mathbf G}
\def\bfH{\mathbf H}
\def\bfI{\mathbf I}
\def\bfJ{\mathbf J}
\def\bfK{\mathbf K}
\def\bfL{\mathbf L}
\def\bfM{\mathbf M}
\def\bfN{\mathbf N}
\def\bfO{\mathbf O}
\def\bfP{\mathbf P}
\def\bfQ{\mathbf Q}
\def\bfR{\mathbf R}
\def\bfS{\mathbf S}
\def\bfT{\mathbf T}
\def\bfU{\mathbf U}
\def\bfV{\mathbf V}
\def\bfW{\mathbf W}
\def\bfX{\mathbf X}
\def\bfY{\mathbf Y}
\def\bfZ{\mathbf Z}

\def\bfw{\mathbf w}

\def\R {{\mathbb R }}
 \def\C {{\mathbb C }}
  \def\Z{{\mathbb Z}}
  \def\H{{\mathbb H}}
\def\K{{\mathbb K}}
\def\N{{\mathbb N}}
\def\Q{{\mathbb Q}}
\def\A{{\mathbb A}}

\def\T{\mathbb T}

\def\bbA{\mathbb A}
\def\bbB{\mathbb B}
\def\bbD{\mathbb D}
\def\bbE{\mathbb E}
\def\bbF{\mathbb F}
\def\bbG{\mathbb G}
\def\bbI{\mathbb I}
\def\bbJ{\mathbb J}
\def\bbL{\mathbb L}
\def\bbM{\mathbb M}
\def\bbN{\mathbb N}
\def\bbO{\mathbb O}
\def\bbP{\mathbb P}
\def\bbQ{\mathbb Q}
\def\bbS{\mathbb S}
\def\bbT{\mathbb T}
\def\bbU{\mathbb U}
\def\bbV{\mathbb V}
\def\bbW{\mathbb W}
\def\bbX{\mathbb X}
\def\bbY{\mathbb Y}

\def\kappa{\varkappa}
\def\epsilon{\varepsilon}
\def\phi{\varphi}
\def\le{\leqslant}
\def\ge{\geqslant}

\begin{center}

{\bf \Large Double cosets for $\SU(2)\times\dots\times\SU(2)$
\\
and outer automorphisms of  free groups}

\vspace{22pt}

{\sc Yuri A. Neretin}%
\footnote{Supported by the grant  FWF, project P19064,
by RosAtom,
and grants
 NWO.047.017.015, JSPS-RFBR-07.01.91209.}

\end{center}

{\small Consider the space of double cosets of
the product of $n$ copies of $\SU(2)$ with respect to the diagonal
subgroup. We get a parametrization
of this space, the radial part of the Haar measure,
and explicit formulas for the actions of the group of outer
automorphisms of the free group $F_{n-1}$
 and of the braid  group of $n-1$ strings.
}

\section{Introduction}

\COUNTERS

{\bf\punct The group $\SU(2)$.}
Denote by $\SU(2)$ the group of unitary
 $2\times 2$-matrices with determinant $=2$.
 A  matrix $g\in\SU(2)$ has the form
$$
g=
\begin{pmatrix}a&b\\ -\ov b& \ov a
  \end{pmatrix},\quad\text{where $|a|^2+|b|^2=1$}
.
$$
Therefore, we can identify the manifold $\SU(2)$ with the 3-dimensional
sphere $S^3$ in $\C^2\subset\R^4$.

\SS


{\bf \punct Double cosets.}
Denote by $\SU(2)$ the group of unitary $2\times 2$
matrices with determinant $=1$. Denote by $G(n)$
the product of $n$ copies of $\SU(2)$.
Elements of $G(n)$ are $n$-tuples
\begin{equation}
(g_1,g_2,\dots,g_n),\qquad \text{where $g_j\in\SU(2)$}
.
\label{eq:tuple}
\end{equation}
Denote
 by $K=K(n)\simeq\SU(2)$
 the diagonal subgroup in $G(n)$;  elements of $K$
 have the form
$$
(h,\dots,h),\qquad \text{where $h\in\SU(2)$.}
$$
The object of the paper is the space  of double cosets
$$
\Pi(n):=K\setminus G/K.
$$
In other words, we consider $n$-tuples (\ref{eq:tuple})
up to the equivalence
$$
(g_1,g_2,\dots,g_n)\sim (h g_1 q,\dots,hg_n q),
\qquad\text{where $h$, $q\in\SU(2)$.}
$$

\smallskip


{\bf\punct Conjugacy classes.}

\begin{observation}
There is a canonical one-to-one correspondence
between $\Pi(n)$ and the space of conjugacy classes of
 $G(n-1)$ with respect to the subgroup $K(n-1)$.
\end{observation}

Indeed,
 $$
 (g_1, g_2,\dots,g_n)\sim (1,g_1^{-1}g_2,\dots, g_1^{-1}g_n)
 $$
  Next,
$$
(1, r_1,\dots,r_{n-1})\sim (1,hr_1h^{-1},\dots,hr_{n-1}h^{-1}).
$$

\smallskip

{\bf\punct Closed polygonal curves on the sphere.}
 Consider the 3-dimensional
sphere $S^3$ endowed with
 the usual angular distance $d(\cdot,\cdot)$.
Fix positive numbers $\theta_1$,\dots, $\theta_{n-1}$.
Consider a closed polygonal curve $A_1A_2\dots A_{n-1}A_1$
in $S^3$ such that
$d(A_{j},A_{j+1})=\theta_j$, $d(A_{n-1},A_1)=\theta_{n-1}$.
Denote by $\cX(\theta)$ the set of all such curves
defined up to proper rotations of the sphere.

\begin{observation}
There is a one-to-one correspondence between
$\cX(\theta)$ and the set of $(n-1)$-tuples
\begin{equation}
(1,r_1,\dots,r_{n-1}),\qquad\text{where $r_j\in\SU(2)$}
\label{eq:rtheta}
\end{equation}
defined up to a simultaneous conjugation
and satisfying the conditions:

\smallskip

--- the eigenvalues of $r_k$ are $e^{\pm i\theta_k}$;

\smallskip

--- $r_1r_2\dots r_{k-1}=1$.
\end{observation}

Indeed, $\SU(2)$ can be considered as a 3-dimensional sphere.
To a tuple (\ref{eq:rtheta}), we assign the polygonal curve
$$
1,\,r_1,\,r_1r_2,\,\dots,\, r_1r_2\dots r_{n-1}=1.
$$

The space $\cX(\theta)$
 (and its analog for $\R^3$ and the Lobachevsky
3-space) became a subject of investigations after Klyachko's
work \cite{Kly}, see e.g. \cite{Kap}, \cite{Jap}.
 Relations of the present work with this literature
is not quite clear for the author. Some other related works are
Fock, Rosly\cite{FR}, Fock \cite{Fock}, Goldman \cite{Gol},
 Dynnikov \cite{Dyn}.

\smallskip


{\bf \punct Spectral forms.}
For a point of $\Pi(n)$, we write the {\it spectral form}
$$
Q(\lambda):=
\det\Bigl(\sum_j \lambda_j g_j\Bigr)
=:\sum_{i,j} s_{ij} \lambda_i\lambda_j
.
$$

We describe  the set $\Xi(n)$ of possible spectral forms.
Namely, they satisfy the conditions:

\smallskip

--- $Q(\lambda)\ge 0$;

\smallskip

--- $\rk Q\le 4$;

\smallskip

--- $s_{jj}=1$.

\smallskip

If $\rk Q=4$, its preimage $\in\Pi(n)$
is a two-point set; we have a branching along the surface
$\rk Q=3$.

Note, that points of the surface $\rk Q=3$
 corresponds to smooth points of
the quotient space $\Pi(n)=K\setminus G/K$; the singular locus
of  $\Pi(n)$ corresponds to the surface $\rk Q\le 2$.

\smallskip


{\bf \punct Radial part of Haar measure.} The group
$G(n)$ is endowed with the Haar measure. We
consider its pushforward to the space $\Xi(n)$.
For $n=3$ we get the usual Lebesgue measure
$ds=ds_{12}\,ds_{13}\,ds_{23}$ on $\Xi(3)$, see \cite{Ner}.
For $n=4$ the measure is given by
$$
\det\bigl(Q(s)\bigr)^{-1/2}\,ds,
$$
where $ds$ is the Lebesgue measure. For $n\ge 5$
the description of the measure  is given in Theorem \ref{th:haar}.

\smallskip


{\bf\punct The group $\Out(F_k)$,}
(see \cite{CV}, \cite{Bes}).
Consider the free group $F_k$ with $k$ generators
$c_1$, \dots, $c_k$.
Denote by $\Aut(F_k)$ the group of
automorphisms of $F_k$. Each automorphism
$\kappa$
is determined by images of the generators:
\begin{equation}
c_j\mapsto \kappa(c_j)=
c_{j_1}^{\epsilon_{j_1}}c_{j_2}^{\epsilon_{j_2}} \dots
\label{eq:action-out}
,
\end{equation}
where $\epsilon_j=\pm1$.  Certainly, the collections
$\bigl\{\kappa(c_j)\bigr\}$ are not arbitrary
(generally, a formula of the type
(\ref{eq:action-out}) determines
 a non-surjective and non-injective map
$F_k\to F_k$).

By the  Nielsen theorem (see \cite{LS}), the group $\Aut(F_k)$
 is generated by the following
  transformations of the set of generators:

\smallskip

a) permutations of generators;

\smallskip

b) the map
$$
c_1\mapsto c_1^{-1},\quad c_2\mapsto c_2,\quad c_3\mapsto c_3,\dots
;
$$

\smallskip

c) the map
$$
c_1\mapsto c_1,\quad c_2\mapsto c_1c_2,\quad
 c_3\mapsto c_3,
\quad
c_4\mapsto c_4,\dots
$$

The group $F_k$ acts on itself by interior automorphisms,
it is a normal subgroup in $\Aut(F_k)$.
We denote by
$$\Out(F_k):=\Aut(F_k)/F_k$$
the {\it group of outer automorphisms of the free group}.

\smallskip


{\bf\punct The action of $\Out(F_{n-1})$ on $\Pi(n)$.}
For a transformation (\ref{eq:action-out}) we write
the following transformation of $\wt\Pi_n$:
$$
\wt r_j=
r_{j_1}^{\epsilon_{j_1}}c_{j_2}^{\epsilon_{j_2}} \dots
$$
In Section 4 we obtain  explicit formulas for
the Nielsen generators.

\smallskip


{\bf\punct Braid groups.} See an introduction in
\cite{Han}, and \cite{Deh}, \cite{DDW}, \cite{Gar}.
 Denote by $\Br_k$ the Artin {\it braid
group}. It has generators $\sigma_1$, \dots $\sigma_{k-1}$
 and relations
\begin{align}
\sigma_j \sigma_{j+1} \sigma_j&=\sigma_{j+1} \sigma_j \sigma_{j+1}
\label{eq:braid-1}
\\
\sigma_i\sigma_j&=\sigma_j\sigma_i\qquad\text{if $|i-j|>1$}
.
\label{eq:braid-2}
\end{align}
There is the following Artin embedding $\Br_k\mapsto\Aut(F_k)$
(see \cite{Han}).
The element $h_j$ corresponds to the transformation
$$
c_j\mapsto c_j c_{j+1} c_j^{-1},\quad c_{j+1}\mapsto  c_j,
$$
other generators are fixed.

There is a characterization (the Artin theorem) of the image of $\Br_k$ in
$\Aut(F_k)$. Namely, $\kappa\in\Out(F_k)$
is  contained
in $\Br_k$ if

\smallskip

1) $\kappa$ sends
each generator $c_j$ to
$$
A_j^{-1} c_{\xi(j)}A_j
$$
where $A_j\in F_k$ and $\xi$ is a permutation of generators.

\smallskip

2) $\kappa$ sends $c_1\dots c_k$ to itself.

\smallskip

In particular, we get the map $\Br_k\to\Out(F_k)$.
It is not injective (see, e.g., \cite{Gar}),
 the kernel
 is generated
by
$$
\Bigl(
(\sigma_1\sigma_2\dots\sigma_{k-1})
(\sigma_1\sigma_2\dots\sigma_{k-2})\dots \sigma_1
\Bigr)^2
.$$

In Section 4 we get explicit
 formulas for the action of generators
of the braid group in the terms of spectral forms.


\smallskip

{\bf\punct Pure braids.}
The relations 
$$
\sigma_j^2=1
$$
together with (\ref{braid-1})--(\ref{braid-2})
determine the symmetric group $S_n$. 
Therefore we get the homomorphism of $\Br_n\to S_n$.
The  kernel is called the {\it group of pure braids}.
The group of pure braids acts on $F_n$ by transformations of the form
$$
c_j\mapsto A_j c_jA_j^{-1}
$$
recall that $c_1\dots c_n\mapsto_1\dots c_n$. This implies the following  observation

\begin{observation}
The group of pure braids acts on the space $\cX(\theta)$ of 
closed polygonal curves.
\end{observation}

\smallskip

{\bf\punct The structure of the paper.}
In Section 2,
we get the characterization  of spectral forms.
 Section 3 contains evaluation of the radial part of the Haar measure.
 In Section 4, we write out actions of the the groups $\Out(F_{n-1})$
 and $\Br_{n-1}$.

 I am grateful to R.~S.~Ismagilov, V.~A.~Fock,
  A.~A.~Rosly, L.~G.~Rybnikov, and I.~A.~Dynnikov
   for discussion of this
 topic.

\section{Spectral forms}

\COUNTERS


{\bf\punct  Spectral forms and the map $\zeta$.}
For any element of $\Pi(n)$
we write out the quadratic form
$$
Q(\lambda_1,\dots,\lambda_n)=\det\bigl(\sum\lambda_j g_j\bigr)=
\det\left(\sum \lambda_j
\begin{pmatrix}a_j& b_j\\-\ov b_j&\ov a_j \end{pmatrix}
\right)
$$
We denote by $\zeta$ the map from $\Pi(n)$ to the space of quadratic
forms.

\begin{proposition}
a) We get a well-defined map from the space $\Pi(n)$
to the space of quadratic forms.

\smallskip

b) Coefficients of $Q$ are real, coefficients in the front
of $\lambda_j^2$ are 1.

\smallskip

c) $Q$ is positive semidefinite.

\smallskip

d) The rank of $Q$ is $\le 4$.
\end{proposition}

{\sc Proof.}
Indeed, for real $\lambda$,
\begin{align*}
Q(\lambda)=\left(\sum \lambda_j a_j\right)
\left(\sum \lambda_j \ov a_j\right)
+
\left(\sum \lambda_j b_j\right)
\left(\sum \lambda_j \ov b_j\right)
=\\=
 \left|\sum \lambda_j a_j\right|^2
+\left|\sum \lambda_j b_j\right|^2
=
\sum \lambda_j^2+ 2 \sum_{i<j} \Re(a_j\ov a_j+b_i\ov b_j)
\end{align*}
and all the statements become obvious.\hfill$\square$

\smallskip


{\bf\punct A description of $\Pi(n)$.}
Denote by $\Xi=\Xi(n)$ the set of all quadratic forms
$Q$ satisfying the conditions of the previous statement.

Obviously,
$$
\zeta(g_1^t,\dots,g_n^t)=\zeta(g_1,\dots,g_n)
,
$$
where $^t$ denotes the transposed matrix.

\begin{theorem}
\label{th:correspondence}
a) The map $\zeta:\Pi(n)\to\Xi(n)$
is surjective.

\smallskip

b) The $\zeta$-preimage of a point $Q\in\Xi$ consists of two
points if $\rk Q=4$ and of one point if
$\rk Q\le 3$.

\smallskip

c) Moreover, $\rk Q\le 3$ iff $(g_1^t,\dots,g_n^t)$ and
$(g_1,\dots,g_n)$ represent one point of $\Pi$.
\end{theorem}

{\sc Proof.}
For a positive semi-definite  quadratic form
$$
Q(\lambda)=\sum_{kl} s_{kl}\lambda_k\lambda_l
$$
 on $\R^n$
there is a collection (a configuration)
of vectors $v_j$ in a Euclidean space such that
$$
\langle v_k,v_j \rangle =s_{kj}
$$
Since $\rk Q\le 4$, this configuration can be realized in $\R^4$.
Since $s_{jj}=1$, these vectors lie on the unit sphere.

Moreover such a configuration $v_j\in\R^4$
is unique up to the action of the orthogonal group
$\OO(4)$.

 Recall that $\SU(2)$ can be identified with
  the 3-dimensional sphere $S^3$,
then $s_{ij}$ are inner products of points of the sphere.
Recall that $\SO(4)\simeq\SU(2)\times\SU(2)/\{\pm 1\}$.
In other words, proper isometries of the sphere
$S^3$ correspond to
the left-right action of $\SU(2)\times\SU(2)$
on $\SU(2)$ (see, e.g. \cite{Zhe}).

If $\rk Q=4$, then $v_j$ are not contained in a
3-dimensional hyperplane. Therefore an improper
isometry of the sphere gives a non-equivalent configuration
in $\Pi(n)$.

If $\rk Q\le 3$, then  the point configuration $v_j$
is contained in a hyperplane. The reflection with
respect to the hyperplane fix this configuration.
\hfill $\square$

\section{The radial part of the Haar measure}

\COUNTERS

{\bf\punct A reduction.}
As we noted above, any element of
the double coset space $\Pi(n)$ can reduced to the form
$(1,g_2,\dots,g_n)$ and $g_j$ are determined up to a simultaneous
conjugation.

\begin{proposition}
Each element of $\Pi(n)$ has a representative
of the form
\begin{multline}
\left\{
\begin{pmatrix} 1&0\\0&1\end{pmatrix},
\begin{pmatrix}e^{i\phi}&0\\0&e^{-i\phi}\end{pmatrix},
\begin{pmatrix}a_2&b_2\\-\ov b_2&\ov a_2 \end{pmatrix},
\dots,
\begin{pmatrix}a_{n-1}&b_{n-1}\\-\ov b_{n-1}&\ov a_{n-1}\end{pmatrix}
\right\},
\\
 \text{where $b_2\ge 0$ and $0\le \phi\le\pi$.}
\label{eq:classes}
\end{multline}
For an element in a general position such a representative is unique.
\end{proposition}

Indeed, after a reduction of $g_2$ to a diagonal
form, we can conjugate our tuple by diagonal matrices.

\smallskip


{\bf\punct The Haar measure on $\SU(2)$.}
We can regard the group $\SU(2)$ as the unit sphere
in the Euclidean space $\C^2$. The Haar measure
on $\SU(2)$ is the usual surface Lebesgue measure
on the sphere. We denote this measure by $dg$.

Recall the following simple facts.

\begin{proposition}
a) The image of the Haar measure under the map
$\begin{pmatrix}a&b\\-\ov b&\ov a \end{pmatrix}\mapsto a$ is
the Lebesgue measure $da\,d\ov a$ on the circle
$|a|\le 1$.

\smallskip

b) Represent $b$ in the form $b=\rho e^{i\theta}$.
Then the image of the Lebesgue measure under the map
$\begin{pmatrix}a&b\\-\ov b&\ov a \end{pmatrix}\mapsto (a,\theta)$
is
$$
\mathrm{const}\cdot d\theta\, da\,d\ov a
.
$$

c) Consider the map taking a matrix $g$ to its
collection of eigenvalues $e^{i\phi}$, $e^{-i\phi}$,
where $0\le \phi\le\pi$.
The image of the Haar measure under the map
$g\mapsto \phi$ is $\sin^2\,\phi d\phi$.

\end{proposition}

\begin{corollary}
The pushforward of the Haar measure in the coordinates
(\ref{eq:classes}) is
$$
\sin^2 \phi\, d\phi\, da_2 \,d\ov a_2\, dg_3\,\dots\,dg_{n}
$$
\end{corollary}


{\bf\punct Coordinates. n-tuples of matrices.}
To be definite, take $n=5$,
\begin{multline*}
\bigl(g_1,g_2,g_3,g_4,g_5\bigr)
:=\\:=
\left( \begin{pmatrix}1&0\\0&1\end{pmatrix},
\begin{pmatrix}e^{i\phi} &0\\0&e^{-i\phi}\end{pmatrix},
\begin{pmatrix}a_1&b_1\\-\ov b_1&\ov a_1\end{pmatrix},
\begin{pmatrix}a_2&b_2\\-\ov b_2 &\ov a_2\end{pmatrix},
\begin{pmatrix}a_3&b_3\\-\ov b_3 &\ov a_3\end{pmatrix}
\right),
\end{multline*}
where $0\le\phi\le\pi$.
We also denote
\begin{align}
a_1&=:x_1+iy_1 \,\, &a_2&=:x_1+i y_2\,\, &a_3&=x_3+iy_3,
\label{eq:phi-et-1}
\\
b_1&=:e^{i\theta_1}\sqrt{1-x_1^2-y_1^2}
\,\,
&b_2&=:e^{i\theta_2}\sqrt{1-x_2^2-y_2^2}
\,\,
&b_3&=:e^{i\theta_3}\sqrt{1-x_3^2-y_3^2}
\label{eq:phi-et-end}
.
\end{align}
For $(g_1,g_2,g_3,g_4,g_5)\in\Pi(n)$ the numbers 
$$
\theta_1-\theta_2, \quad \theta_2-\theta_3, \quad \theta_1-\theta_3
$$
 make sense
(but not $\theta_1$, $\theta_2$, $\theta_3$ themselves).

\smallskip


{\bf\punct Coordinates. Spectral forms.}
Consider the  spectral form and denote
its coefficients in the following way:
\begin{multline*}
\det
\left(\sum \lambda g_1+ \mu g_2+\nu_1 g_3+\nu_2 g_4+\nu_3 g_5
 \right)
=:
\\
=:
\lambda^2+\mu^2+\nu^2+ 2p\lambda \mu\,\,
+2 q_1\lambda\nu_1+2q_2\lambda \nu_2 + 2q_3\lambda \nu_3
\\+2r_1\mu\nu_1+2r_2\mu\nu_2 +2r_3\mu\nu_3\,\,
+2t_{12}\nu_1\nu_2+2t_{13}\nu_1\nu_3+2t_{23}\nu_2\nu_3
.
\end{multline*}
The matrix of the form is
\begin{equation}
\Delta=
\begin{pmatrix}
1&p&q_1&q_2&q_3\\
p&1&r_1&r_2&r_3\\
q_1&r_1&1 &t_{12}&t_{13}\\
q_2&r_2&t_{12}&1&t_{23}\\
q_3&r_3&t_{13}&t_{23}&1
\end{pmatrix}
\label{eq:Delta}
.
\end{equation}

Then
\begin{align}
\label{eq:pqrt-1}
p&=\cos\phi
\\
q_j&=x_j
\\
r_j&= x_j\cos\phi+y_j\sin\phi
\\
t_{ij}&=
 x_ix_j+y_iy_j+\sqrt{1-x_i^2-y_i^2}\sqrt{1-x_j^2-y_j^2}
 \cos(\theta_i-\theta_j)
 \label{pqrt-2}
\end{align}
It is easy to write the inverse map:
\begin{align}
\label{eq:arccosp}
\phi&=\arccos p
\\
x_j&=q_j
\\
y_j&=\frac{r_j-q_jp} {\sqrt{1-p^2}}
\\
\theta_i-\theta_j&=
\pm\arccos \frac{\det
\begin{pmatrix} 1&p&q_i\\p&1&r_i\\q_j&r_j&t\end{pmatrix}}
     {\det^{1/2}
     \begin{pmatrix}1&p&q_i\\p&1&r_i\\q_i&r_i&1\end{pmatrix}
       \det^{1/2}\begin{pmatrix}1&p&q_j\\p&1&r_j\\q_j&r_j&1\end{pmatrix}}
\label{eq:costheta}
\end{align}

The last formula requires some calculations. For this reason,
we present some intermediate formulas:
\begin{align}
1-x_1^2-y_1^2=
\frac{\det\begin{pmatrix}1&p&q_1\\p&1&r_1\\q_1&r_1&1\end{pmatrix}}
{1-p^2}
,
\label{eq:1x2y2}
\\
t-x_1x_2-y_1y_2=\frac
{\det\begin{pmatrix} 1&p&q_1\\p&1&r_1\\q_2&r_2&t\end{pmatrix}}
{1-p^2}.
\nonumber
\end{align}

Note, that we can not reconstruct the sign of $\theta_i-\theta_j$
from the formula (\ref{eq:costheta}). Recall that
the substitution
$$
\theta_1\mapsto-\theta_1,
\quad \theta_2\mapsto-\theta_2,
\quad
\theta_3\mapsto-\theta_3
$$
corresponds to the simultaneous transposition
$$
(g_1,g_2,g_3,g_4,g_5)\mapsto (g_1^t,g_2^t,g_3^t,g_4^t,g_5^t)
.
$$


{\bf\punct What happens if we forget $t_{23}$?}
Let we know $p$, all $q_j$, all $r_j$, and $t_{12}$, $t_{13}$.
Then we can reconstruct $\phi$, $x_j$, $y_j$ and
$$
\cos(\theta_1-\theta_2),\quad \cos(\theta_1-\theta_3)
.
$$
Without loss of a generality, we can assume $\theta_1=0$.
Then we know  $\pm\theta_2$ and $\pm\theta_3$ and there are
two possible variants for $|\theta_2-\theta_3|$.

It can be readily checked that there exist $h\in\SU(2)$ such that
$$
h^{-1}g_2 h=g_2^t,
\quad h^{-1}g_3 h=g_3^t
$$
(recall that $g_1$ is the unit matrix). Then
without $t_{23}$ we can not distinguish
\begin{equation}
(g_1,g_2,g_3,g_4,g_5)\quad \text{and}
\quad
(g_1,g_2,g_3,g_4,h g_5^t h^{-1}).
\label{eq:ggggg-ggggg}
\end{equation}


{\bf\punct The radial part of the Haar measure.
The cases $n=3$, $n=4$.}

\begin{theorem}
\label{th:3-4}
a) Let $n=3$. The pushforward of the Haar measure
under the map $\zeta:\Pi(3)\to\Xi(3)$ is%
\footnote{See, also \cite{Gol}.}
$$
\mathrm{const}\cdot
dp\,dq_1\,dr_1
.
$$

b) Let $n=4$. Then the image of the Haar measure under the map
$\zeta:\Pi(4)\to\Xi(4)$ is
$$
\mathrm{const}\cdot
 \det \begin{pmatrix}
1&p&q_1&q_2\\
p&1&r_1&r_2\\
q_1&r_1&1 &t_{12}\\
q_2&r_2&t_{12}&1
\end{pmatrix}^{-1/2}
  \,dp\,dq_1\,dq_2\,dr_1\,dr_2\,dt_{12}
  .
$$
\end{theorem}

{\sc Proof.} Consider the case $n=4$. The radial part of
the Haar measure in the coordinates $\phi$, $x_1$, $y_1$,
$x_2$, $y_2$, $\theta$ is given by
$$
\mathrm{const}\cdot\sin^2\phi \,d\phi\,dx_1\,dy_1\,dx_2\,dy_2\,d\theta
$$
Next, we must write the Jacobian of the map
(\ref{eq:arccosp})--(\ref{eq:costheta}).
Evidently, the Jacobian is
$$
\frac{\partial \phi}{\partial p}
\cdot
\frac{\partial \theta}{\partial t}
,$$
this can be easily evaluated.

\smallskip

{\bf\punct The Haar measure, general case.}
For an $n$-tuple $(g_1,\dots,g_n)$ consider
its spectral form
$$
\det\Bigl(\sum_j \lambda_jg_j\Bigl)
=:\sum_j \lambda_j^2+2\sum_{i<j}s_{ij}\lambda_i\lambda_j
$$

\begin{theorem}
\label{th:haar}
a) The coefficients $s_{12}$, $s_{13}$, $s_{23}$ are
distributed as
$$
\mathrm{const}\cdot
ds_{12}\,ds_{13}\,ds_{23}
$$
in the domain
$\begin{pmatrix}
1&s_{12}&s_{13} \\
s_{12}&1&s_{23}\\
s_{13}&s_{23}&1\end{pmatrix}\ge0$.

\smallskip

b) For fixed $s_{12}$, $s_{13}$, $s_{23}$ in a general position,
a vector $v_j:=\begin{pmatrix} s_{1j}&s_{2j}&s_{3j}\end{pmatrix}$
is distributed as
$$
\mathrm{const}\cdot
\det(\Delta_j)^{-1/2} ds_{j1}\,ds_{j2}\,ds_{j3}
,$$
where
$$
\Delta_j=\begin{pmatrix}
1&s_{12}&s_{13}&s_{1j} \\
s_{12}&1&s_{23}&s_{2j}\\
s_{13}&s_{23}&1&s_{3j}\\
s_{1j}&s_{2j}&s_{3j}&1
\end{pmatrix}
.$$
A vector $\begin{pmatrix} s_{1j}&s_{2j}&s_{3j}\end{pmatrix}$
ranges in the domain $\Delta_j\ge 0$.

\smallskip

c) The random variables $v_4$, $v_5$, \dots, $v_n$ are independent.

\smallskip

d) Let us fix $s_{1j}$, $s_{2j}$, $s_{3j}$ for all $j$.
For such a collection in a general position there are 2
eqiuprobable
variants of a choice of $s_{4j}$. These samplings are independent
for $j=5$, $6$, \dots.

\smallskip

e) Let us fix $s_{1j}$, $s_{2j}$, $s_{3j}$ for all $j$ and
fix
$s_{4j}$. Then a.s. all other variables $s_{ij}$ are uniquely
determined.
\end{theorem}

{\sc Proof.} The statements a)-b) are a rephrasing of
Theorem (\ref{th:3-4}).

Next, for fixed $g_1$, $g_2$, $g_3$ the matrices
('random variables') $g_3$, $g_4$, \dots are independent.
A matrix $g_j$ determines a vector $v_j$, and $g_j$ is uniquely
determined by a vector $v_j$. This proves c).

Now we write the matrix of a spectral form
$$
\begin{pmatrix}
1&s_{12}&s_{13}&s_{14}&s_{15}&s_{16}&\dots
\\
s_{12}&1&s_{23}&s_{24}&s_{25}&s_{26}&\dots
\\
s_{13}&s_{23}&1 &s_{34}&s_{35}&s_{36}&\dots
\\
s_{14}&s_{24}&s_{34}&1&?_1& ?_2 &\dots
\\
s_{15}&s_{25}&s_{35}&?_1&1& *&\dots
\\
s_{16}&s_{26}&s_{36}&?_2 &* &1&\dots
\\
\vdots&\vdots&\vdots&\vdots&\vdots&\vdots&\ddots
\end{pmatrix}
.
$$
The left upper $5\times 5$ minor is zero. This gives
a quadratic equation for $s_{45}$. Both the solutions are
admissible, they correspond to collections of matrices
given by (\ref{eq:ggggg-ggggg}).
We repeat the same operation for $s_{46}$ etc.

In the terminology of Subsections 3.4--3.5,
the knowledge of the first three rows gives us
  $\pm \theta_1$, $\pm\theta_2$, $\pm\theta_3$ etc.
  All these choices are equiprobable.

Having all $s_{4j}$,
we get a unique way to complete the matrix to the matrix
of rank 4.

\section{Actions of $\Out(F_{n-1})$ and braid group}

\COUNTERS

\smallskip


{\bf\punct The action of $\Out(F_{n-1})$ on $\Pi(n)$.}
We regard $\Pi(n)$ as the set of collections
$(1, g_2, \dots,g_n)$ defined up to a simultaneous  conjugation.

The Nielsen transformations act on $\Pi[n)$ in the obvious way.
We get:

\smallskip

a) permutations of $g_j$;

\smallskip

b) the transformation $g_2\mapsto g_2^{-1}$;

\smallskip

c) the map $(g_2,g_3, g_4,\dots)\to (g_2, g_2g_3, g_4,\dots)$.

\smallskip

To be definite, take $n=5$. The action of permutations is obvious.

\begin{proposition}
The transformation $g_2\mapsto g_2^{-1}$
corresponds to the map $\Xi(5)\to\Xi(5)$ given by
\begin{align*}
\wt p&= p;
\\
\wt q_j&=q_j;
\\
\wt t_{ij}&=t_{ij}
\\
\wt r_j&=-r_j+2pq_j
\end{align*}
\end{proposition}

{\sc Proof.}
In the notation (\ref{eq:phi-et-1})--(\ref{eq:phi-et-end}),
(\ref{eq:pqrt-1})--(\ref{pqrt-2}),
we have
$$
\wt r_1=x_1\cos\phi-x_1\sin\phi=r_1-2y_1\sin\phi
.
$$
On the other hand,
$$
y_1\sin\phi=r_1-x_1\cos\phi=r_1-pq_1
.
$$
This completes the calculation. \hfill $\square$

\begin{theorem}
\label{th:nielsen-2}
The transformation
$$
(1,g_2,g_3,g_4,g_5)\mapsto (1,g_2,g_2g_3,g_4,g_5)
$$
corresponds to the map $\wt\Xi(5)\to\wt\Xi(5)$ given by
\begin{align}
\wt p&=p;
\nonumber
\\
\wt q_2&=q_2,\quad \wt q_3=q_3;
\nonumber
\\
\wt q_1&=-r_1+2pq_1;
\nonumber
\\
\wt r_2&=r_2,\quad\wt r_3=q_3;
\nonumber
\\
\wt t_{1j}&=p\, t_{1j}-q_j r_1+q_1 r_j
\mp\det
\begin{pmatrix}
1&p&q_1&q_j\\
p&1&r_1&r_j\\
q_1&r_1&1&t_{1j}\\
q_j&r_j&t_{1j}&1
\end{pmatrix}^{1/2}\!\!\!\!\!\!,\,\,\,\text{where $j=2$, $3$}.
\label{eq:first-long}
\\
\wt t_{23}&=t_{23}
\nonumber
\end{align}
\end{theorem}

The group $\Out(F_n)$ acts on $\wt\Xi(5)$ and not on
$\Xi(5)$ and the choice of signs $\mp$ requires explanations.
They are given below.

\smallskip

{\sc Proof.} First, we write the coefficients of the spectral
form for the transformed collection. Only the variables
$q_1$, $r_1$, $t_{12}$, $t_{13}$ change.
We have
\begin{multline*}
\wt r_1=\wt x_1\cos\phi+\wt y_1\sin\phi
=\\=
(x_1\cos\phi-y_1\sin\phi)\cos\phi+
     (y_1\cos\phi+x_1\sin\phi)\sin\phi=x_1=q_1
.
\end{multline*}
and
\begin{equation}
\wt q_1=\wt x_1=x_1\cos\phi-y_1\sin\phi=
q_1p-(r_1-q_1p)=-r_1+2q_1p
\label{eq:wtq1}
.
\end{equation}
The evaluation of $\wt t_{1j}$ is heavier,
$$
\wt t_{12}=\wt x_1 \wt x_2+\wt y_1 \wt y_2 +
\sqrt{1-\wt x_1^2-\wt y_1^2}\sqrt{1-\wt x_2^2-\wt y_2^2}
\cos(\wt\theta_1-\wt\theta_2)
.
$$
By construction, $\wt x_2=x_2$, $\wt y_2=y_2$,
$\wt \theta_2=\theta_2$. Next,
$$
\wt a_1=\wt x_1+i\wt y_1=(x_1+iy_1)\,e^{i\phi}
$$
and therefore
$$
1-\wt x_1^2-\wt y_1^2=1-x_1^2-y_1^2
.
$$
Also, $\wt\theta_1=\theta_1+\phi$.
Denote $\theta:=\theta_1-\theta_2$.
Thus,
$$
\wt t_{12}=\wt x_1 x_2+\wt y_1  y_2 +
\sqrt{1- x_1^2- y_1^2}\,\sqrt{1- x_2^2- y_2^2}
\bigl(
\cos\theta\cos\phi-\sin\theta\sin\phi
\bigr)
.
$$
The variable $\wt x_1$ is evaluated in (\ref{eq:wtq1}),
$$
\wt y_1=x_1\sin\phi+y_1\cos\phi=
\frac{r_1-q_1p}{\sqrt{1-p^2}}\cdot p
+
q\sqrt{1-p^2}
.
$$
We use formula (\ref{eq:1x2y2}) for
square roots and formula (\ref{eq:costheta}) for $\cos\theta$.
After this, we can evaluate $\sin\theta$,
$$
\sin^2\theta
=
\frac{(1-p^2)\cdot \det
\begin{pmatrix}
1&p&q_1&q_2\\
p&1&r_1&r_2\\
q_1&r_1&1&t_{12}\\
q_2&r_2&t_{12}&1
\end{pmatrix}}
{\det\begin{pmatrix}
1&p&q_1\\
p&1&r_1\\
q_1&r_1&1
\end{pmatrix}
\det\begin{pmatrix}
1&p&q_2\\
p&1&r_2\\
q_2&r_2&1
\end{pmatrix}
}
.$$
After this, we get a unexpectedly long chain of cancelations
and get the desired formula.
\hfill $\square$.

\smallskip

{\sc Choice of signs.} We use formula (\ref{eq:costheta})
and find
$$
\pm(\theta_1-\theta_2),\quad
\pm(\theta_1-\theta_3),\quad
\pm(\theta_2-\theta_3)
$$
These numbers must be consistent, in fact only two variants
are possible (this corresponds to a choice of a sheet
of the covering map $\wt\Xi\to\Xi$). Now let we have chosen
the signs. Then we take 'minus' in
(\ref{eq:first-long}) if $(\theta_1-\theta_j)\ge 0$.
Otherwise, we take 'plus'.

\smallskip


{\bf\punct The action of the braid group.}

\begin{lemma}The transformation
$$
(1,g_2,g_3,g_4,g_5)\mapsto
(1,g_2g_3g_2^{-1},g_2,g_4,g_5)
$$
corresponds to the map $\wt\Xi(5)\to\wt\Xi(5)$
given by
\begin{align*}
\wt p&= p
\\
\wt q_k&=q_k, \quad\text{where $k=1$, $2$, $3$;}
\\
\wt r_k&=r_k, \quad\text{where $k=1$, $2$, $3$;}
\\
\wt
t_{1j}&=
t_{1j}-2\cdot
\det\begin{pmatrix}
1&p&q_1\\
p&1&r_1\\
q_j&r_j&t_{1j}
\end{pmatrix}
\mp 2p\cdot
\det
\begin{pmatrix}
1&p&q_1&q_j\\
p&1&r_1&r_j\\
q_1&r_1&1&t_{1j}\\
q_j&r_j&t_{1j}&1
\end{pmatrix}^{1/2},
\\
\wt t_{23}&=t_{23}.
\end{align*}

\end{lemma}

{\sc Proof.} We  evaluate
\begin{multline*}
\wt t_{12}=
\Re(a_1\ov a_2+b_1\ov b_2 e^{2i\phi})
=\\=
x_1x_2+y_1y_2+
\sqrt{1-x_1^2-y_1^2}\sqrt{1-x_2^2-y_2^2}
(\cos\theta\cos2\phi-\sin\theta\sin 2\phi)
\end{multline*}
as in the previous proof.\hfill$\square$

\smallskip

Now we can write the action of the braid group.
To write  formulas for generators, we represent
 the matrix of the spectral form  as
$$
\begin{pmatrix}
1&p_1&p_2&p_3&\dots\\
p_1&1&h_{12}&h_{13}&\dots\\
p_2&h_{12}&1&h_{23}&\dots\\
p_3&h_{13}&h_{23}&1&\dots\\
\vdots&\vdots&\vdots&\vdots&\ddots
\end{pmatrix}
$$
We also set $h_{ij}:=h_{ji}$.
Then the formula for a generator
$\sigma_k$ of the braid group is
\begin{align*}
\wt h_{ij}&=h_{ij},\quad\text{if $i$, $j\ne k$, $k+1$}
\\
\wt p_i&=p_{i},\quad\text{if $i\ne k$, $k+1$}
\\
\wt p_k&=p_{k+1}
\\
\wt p_{k+1}&=p_k
\\
\wt h_{(k+1)j}&=h_{kj},
\end{align*}
and
\begin{multline}
\wt
h_{kj}=h_{(k+1)j}-2\cdot\det
  \begin{pmatrix}
   1&p_k& p_{k+1}\\
   p_k&1&h_{k(k+1)}\\
   p_{j}&h_{kj} &h_{(k+1)j}
   \end{pmatrix}
   -\\
   -
   2p_k\cdot
   \det
    \begin{pmatrix}
   1&p_k& p_{k+1}&p_j\\
   p_k&1&h_{k(k+1)}&h_{kj}\\
   p_{k+1}&h_{k(k+1)}&1&h_{(k+1)j}
   \\
   p_{j}&h_{kj} &h_{(k+1)j}&1
   \end{pmatrix}^{1/2}
\end{multline}

{\tt Math.Dept., University of Vienna,

 Nordbergstrasse, 15,
Vienna, Austria

\&

Institute for Theoretical and Experimental Physics,

Bolshaya Cheremushkinskaya, 25, Moscow 117259,
Russia

\&

Mech.Math. Dept., Moscow State University,
Vorob'evy Gory, Moscow

e-mail: neretin(at) mccme.ru

URL:www.mat.univie.ac.at/$\sim$neretin

wwwth.itep.ru/$\sim$neretin
}

\end{document}